\documentclass{amsart}

\newtheorem{theorem}{Theorem}[section]
\newtheorem{lemma}[theorem]{Lemma}

\theoremstyle{definition}
\newtheorem{definition}[theorem]{Definition}
\newtheorem{example}[theorem]{Example}

\theoremstyle{remark}
\newtheorem{remark}[theorem]{Remark}

\numberwithin{equation}{section}

\usepackage{amsthm,amsmath,amssymb,xcolor,stmaryrd}
\usepackage[inline]{enumitem}
\usepackage{multicol,nicefrac}

\newtheorem{corollary}[theorem]{Corollary}
\newtheorem{proposition}[theorem]{Proposition}

\makeatletter
\newcommand*{\ineq}[2][]{%
  \begingroup
    \refstepcounter{equation}%
    \ifx\\#1\\%
    \else
      \label{#1}%
    \fi
    \relpenalty=10000 %
    \binoppenalty=10000 %
    \@eqnnum \ \ensuremath{%
      \hskip 20pt #2%
    }%
  \endgroup
}
\makeatother

\newcommand{\vnf}{\nfpar\omega}
\newcommand{\penum}[1]{\check #1}
\newcommand{\penu}[1]{ #1_{-1}}
\newcommand{\xcases}[2]{\nicefrac #1  #2}
\newcommand{\fs}[2]{[#2]#1}
\newcommand{\Fix}[1]{{\rm Fix}(#1)}
\newcommand{\atr}{{\sf ATR}_0}
\newcommand{\fsi}[2]{\lb #2\rb #1}
\newcommand{\ok}[2]{o_{#2}#1}
\newcommand{\peq}{\preccurlyeq}
\newcommand{\seq}{\succcurlyeq}
\newcommand{\good}[2]{{\mathbb G}_{#2}#1}
\newcommand{\lb}{\llbracket}
\newcommand{\rb}{\rrbracket}

\newcommand{\bch}[2]{ \langle #2 \rangle #1 }
\newcommand{\bases}[2]{#2/#1}
\newcommand{\david}[1]{}

\newcommand{\putaway}[1]{}

\newcommand{\cnfsym}{\approx}

\newcommand{\cnf}{\cnfsym_k}
\newcommand{\gnfsym}{\cong}
\newcommand{\gnf}{\gnfsym_k}
\newcommand{\gknfsym}{\equiv}

\newcommand{\nfpar}[1]{\gknfsym_{#1}}
\newcommand{\enfpar}[1]{\gnfsym_{#1}}
\newcommand{\gknf}{\nfpar{k} }

\newcommand{\pk}{\langle {\bases k \omega} \rangle}
\newcommand{\ps}[1]{\langle {\bases {#1}{\omega}} \rangle }

\newcommand{\pskpe}{\ps{k+1}}

\newcommand{\ot}{{\phi}}

\newcommand{\al}{{\alpha}}

\newcommand{\be}{{\beta}}

\newcommand{\ga}{\gamma}

\newcommand{\la}{\lambda}

\newcommand{\om}{\omega}

\newcommand{\longversion}[1]{}
\newcommand{\shortversion}[1]{
#1
}

\begin{document}



\address{}
\curraddr{}
\email{}
\thanks{}

\author{}
\address{}
\curraddr{}
\email{}
\thanks{}


\date{}

\dedicatory{}

\commby{}

\begin{abstract}
The classical Goodstein process gives rise to long but finite sequences of natural numbers whose termination is not provable in Peano arithmetic.
In this manuscript we consider a variant based on the Ackermann function.
We show that Ackermannian Goodstein sequences eventually terminate, but this fact is not provable using predicative means.
\end{abstract}

\title[The Ackermannian Goodstein process]{
Predicatively unprovable termination of the Ackermannian Goodstein process\\
}

\author[T.~Arai]{Toshiyasu Arai}
\author[D.~Fern\'andez]{David Fern\'andez-Duque}
\author[S.~Wainer]{Stanley Wainer}
\author[A.~Weiermann]{Andreas Weiermann}

\maketitle

\section{Introduction}

Among the greatest accomplishments of mathematical logic in the first half of the twentieth century was the identification of true arithmetical statements unprovable in Peano arithmetic ($\sf PA$): the consistency of $\sf PA$, due to G\"odel \cite{Godel1931}, and transfinite induction up to the ordinal $\varepsilon_0$, due to Gentzen \cite{Gentzen1936}.
However, such statements do not clarify whether incompleteness phenomena should be pervasive in other disciplines such as combinatorics or number theory.

In contrast, Goodstein's principle \cite{Goodsteinb} is a purely number-theoretic statement simple enough to be understood by a high school student yet unprovable in $\sf PA$.
While the statement makes no reference to ordinals, Goodstein's proof uses the well-foundedness of $\varepsilon_0$ \cite{Goodstein1944}.
Kirby and Paris' later independence proof shows that this use of $\varepsilon_0$ is, in some sense, essential \cite{Kirby}.
These developments paved the road for the discovery of other combinatorial statements independent from $\sf PA$, including the Paris-Harrington theorem \cite{ParisHarrington} and other Ramsey-style principles \cite{Pelupessy,vanHoof}.

A modern presentation of Goodstein's result and proof consists of three ingredients.
First, we need a notion of {\em normal form} for representing natural numbers.
Say that a natural number $m>0$ has {\em base-$k$ exponential normal form $k^a+b$} if $m=k^a+b$ and $k^a \leq m < k^{a+1} $.
By iteratively applying this definition to $b$ we obtain a standard base-$k$ representation of $m$.
Note that the exponents themselves can be recursively written in base $k$: for example, $21 = 2^{2^2} + 2^2 + 2^0$ in base $2$.

The second ingredient is the base-change operation, which replaces every occurrence of $k$ by some $\ell>k$: formally, $\bch 0{\bases k\ell} = 0$ and $\bch {(k^a + b)}{\bases k\ell} = \ell^{\bch a{\bases k\ell}} + \bch b{\bases k\ell}$. For brevity we write $\bch {m}{k+1}$ instead of $\bch m{\bases k{k+1}}$.
With this, given a natural number $m$, we define a sequence $({\rm G}_i m)_{m<\alpha}$, where $\alpha \leq \omega$, recursively by
${\rm G}_0m = m$, ${\rm G}_{i+1}m = \bch {{\rm G}_im}{i+3} - 1 $ if ${\rm G}_im > 0$, and setting $\alpha = i+1 $ if ${\rm G}_im = 0$; if no such $i$ exists, $\alpha = \omega$.
Goodstein's principle is then the following.

\begin{theorem}[Goodstein]\label{theoGood}
For every $m\in \mathbb N$ there is $i\in\mathbb N$ such that ${\rm G}_im = 0$.
\end{theorem}

For the proof we need an additional ingredient: an interpretation of the Goodstein process in terms of ordinals below $\varepsilon_0$.
Every non-zero $\xi<\varepsilon _ 0$ can be written uniquely in the form $\omega^\alpha + \beta$ with $\alpha,\beta<\xi$; this is the {\em Cantor normal form} of $\xi$.
Given natural numbers $m,k$ with $m=k^a+b$ in base-$k$ exponential normal form and $k>1$, we may recursively define an ordinal $\pk m = \omega^{\pk a} + \pk b$, setting $\pk 0 = 0$.
Induction on $m$ shows that $\pk m<\varepsilon_0$ is well-defined.
Now, let $({\rm G}_im)_{i<\alpha}$ be the Goodstein process for $m$, and for $i<\alpha$ define $\ok mi = \ps {i+2}{\rm G}_im $.
It is not too hard to check that $ \ok m0 > \ok m1 > \ok m2 >\ldots . $
By the well-foundedness of $\varepsilon_0$ there can be no such infinite sequence, hence $\alpha<\omega$, as needed.

By coding finite Goodstein sequences as natural numbers in a standard way, Goodstein's principle can be formalized in the language of arithmetic.
However, this formalized statement is unprovable in $\sf PA$. \shortversion{This can be proven by showing that the Goodstein process takes at least as long as stepping down the {\em fundamental sequences} below $\varepsilon_0$; these are canonical sequences\footnote{We write $\fs \xi n$ instead of the more standard $\xi[n]$ since it will ease notation later.} $(\fs \xi n)_{n<\omega}$ such that $\fs \xi n< \xi$ and for limit $\xi$, $\fs \xi n \to \xi$ as $n\to \infty$.
For standard fundamental sequences below $\varepsilon_0$, $\sf PA$ does not prove that the sequence $\xi > \fs \xi 2 > \fs{\fs \xi 2} 3 > \fs{  \fs{\fs \xi 2} 3 } 4 \ldots $ is finite.\medskip}

Exponential notation falls short when attempting to represent large numbers that arise in combinatorics such as Graham's number \cite{graham1969}.
These numbers may instead be written in terms of fast-growing functions such as the Ackermann function $A_a(k,b)$; see Definition \ref{defAckermann}.
Here $k$ is regarded as the `base' and $a,b$ as parameters.
There will typically be several ways to write a number in the form $A_a (k,b) +c$, so a suitable normal form is chosen in Definition \ref{defNF}.
Our normal forms are based on iteratively approximating $m$ via a `sandwiching' procedure.
With this we can define the {\em Ackermannian Goodstein process.} This process always terminates, although the proof now uses the well-foundedness of the Feferman-Sch\"utte ordinal $\Gamma_0$, regarded by Feferman \cite{Feferman} as representing the limit of {\em predicative mathematics} \cite{Weyl}, which denies the existence of the real line as a completed totality.
The ordinal $\Gamma_0$ is the proof-theoretic ordinal of the theory $\atr$ of {\em arithmetical transfinite recursion.}

Let us give a brief description of $\atr$; a more formal treatment can be found in \cite{Simpson:2009:SubsystemsOfSecondOrderArithmetic}.
\shortversion{Second-order arithmetic adds set-variables $X,Y,Z,\ldots$ and predicates $t \in X$ to the language of $\sf PA$. Quantifiers may range over first- or second-order variables.
The theory ${\sf ACA}_0$ of arithmetical comprehension includes axioms for arithmetical operations along with the induction axiom
and the {\em comprehension} scheme stating that $\{x\in \mathbb N: \Phi(x)\}$ is a set, where $\Phi$ is {\em arithmetical,} i.e.~it contains no second-order quantifiers.
Given a linear order $\prec$ and a set $X$ of ordered pairs $\langle\lambda,n\rangle$, let $X_{\lambda} = \{ n\in \mathbb N : \langle \lambda,n\rangle\in X\}$, and $X_{{\prec} \lambda} = \{ \langle \xi,n\rangle \in X: \xi \prec \lambda\}$. The theory $\atr$ is ${\sf ACA}_0$ plus the scheme stating that if $\prec$ is a well-order, there exists a set $X$ such that $\forall n,\lambda  \big (n\in X_\lambda\leftrightarrow \Phi(n,X_{{\prec}\lambda}) \big )$, with $\Phi$ arithmetical. $\atr$ is conservative over systems of predicative mathematics \cite{Feferman:1964:SystemsOfPredicativeAnalysis,
Simpson:1985:FriedmansResearh}.
}

\section{Ackermannian Goodstein Sequences}

In this section we establish the basic definitions needed to state our main results.
\longversion{
All notions will be revisited and elaborated on in later sections.
We begin by defining the version of the Ackermann function we will be working with.}

\begin{definition}\label{defAckermann}
Given $a,b,k \in \mathbb N$ with $k\geq 2$, we define $A_a (k,b) \in \mathbb N$ as follows.
Fix $k$ and let us write $A_a b$ instead of $A_a (k,b)$.
Define as an auxiliary value $A_a (-1) = 1$.
Then, $A_a b$ is given recursively by
$A_0 b = k^b $ and
$A_{a+1} b = A_a^k A_{a+1} (b-1)$.
\end{definition}

Note that $A_0(k,0) = 1$ regardless of $k$.
Aside from some trivial cases, the Ackermann function is strictly increasing on all parameters, as can be verified by an easy induction; we leave the details to the reader.

\begin{lemma}\label{lemmBoundMajorize}
Let $a \leq a'$, $b \leq b'$, and $2\leq k\leq k'$ be natural numbers. Then,
\[\max\{a,b\}< A_a (k,b) \leq A_{a'} (k',b').\]
If moreover $a+b+k<a'+b'+k'$ and $a' + b' > 0$, then $A_a (k,b) < A_{a'} (k',b')$.
\end{lemma}

\longversion{
Next we define normal forms based on the Ackermann function.
Fix a `base' $k$, and write $A_ab$ instead of $A_a(k,b)$.
The general idea is to represent $m>0$ canonically in the form $A_ab + c$.
It is tempting to choose $a$ maximal so that there is $ b$ with $A_a b \leq m < A_a(b+1)$.
However, in this case $c$ may still be quite large; so large, in fact, that there is $a'<a$ with $A_{a'}A_a b \leq m$.
In this case, $A_{a'}A_a b $ is a better approximation to $m$ than $A_ab$; and, indeed, we may choose $b'$ maximal so that $A_{a'}b'\leq m$.
We then have that
\[A_ab < A_{a'} b' \leq m < A_{a'} (b'+1) \leq A_a (b+1),\]
`sandwiching' $m$ between better and better approximations.
Continuing in this fashion, we can find the `best' approximation to $m$; this will be the basis for our normal forms.
}

\shortversion{
In order to canonically represent $m$ in the form $A_ab+c$, we first choose $a_1$ maximal such that $A_{a_1}0\leq m $, then $b_1$ maximal so that $A_{a_1}b_1 \leq m$.
However, $A_{a_1}b_1 $ might be much smaller than $m$, even small enough that there is $a_2$ with $A_{a_2}A_{a_1}b_1 \leq m$.
In this case, $A_{a_2}b_2$ is a `better' approximation to $m$, where $b_2 \geq A_{a_1}b_1$ is maximal.
Continuing this process, we reach the `best' approximation $A_{a_n}b_n \leq m$.
}

\begin{definition}\label{defNF}
Fix $k\geq 2$ and let $A_x y = A_x (k,y)$.
Given $m , a, b, c \in \mathbb N$ with $m>0$, we define $A_a b + c $ to be the {\em $k$-normal form} of $m$, in symbols $m \equiv _k A_a b + c$, if $m =  A_a b + c$ and there exist sequences $a_1 , \ldots, a_n$ of {\em sandwiching indices,}, $b_1,\ldots, b_n$ of {\em sandwiching arguments} and $m_0,\ldots m_n$ of {\em sandwiching values} such that for $i<n$,
\begin{multicols}2

\begin{enumerate}[wide]

\item $m_0 = 0$;

\item $ A_{a_{i+1}} m_i \leq m <  A_{a_{i+1} + 1} m_i$;

\item $ A_{a_{i+1}}  b_{i+1} \leq m <  A_{a_{i+1}} (b_{i+1} + 1)$;

\item $m_{i+1} = A_{a_{i+1}} {b_{i+1}}$;

\item $A_0 m_n > m$, and

\item $a = a_n$ and $b = b_n$.

\end{enumerate}
\end{multicols}
\noindent We denote the sequence of pairs $(a_i,b_i)$ by
$ ( A_{a_i}b_i)_{i=1}^n $
and call it the {\em $k$-sandwiching sequence} of $m$.
If $m = 0$, $m$ has empty sandwiching sequence and $k$-normal form $0$.
\end{definition}

We write simply {\em sandwiching sequence} when $k$ is clear from context.
For our purposes $m \gknf A_ab+ c$ is viewed as a relation between the numbers $a,b,c,m,k$; a different approach where $A_ab+c$ is regarded as a formal term is briefly discussed in Appendix \ref{secConc}.
Every positive integer has a unique $k$-sandwiching sequence and hence a unique normal form. This is because $a_{i+1}$ and $b_{i+1}$ are unique when defined since $A_xy$ is strictly increasing on both variables, and $(m_i)_{i\leq n}$ is strictly increasing (see Lemma \ref{lemmSandwichingMonotone}). Thus we must have $n\leq m $ and $m$ indeed has a finite sandwiching sequence.
It is not hard to see that $1 \gknf A_00$ with sandwiching sequence $(A_00)$.

\begin{example}\label{exHard}
Let us compute the $2$-normal form and the sandwiching sequence $(A_{a_i}b_i)_{i=1}^n$ for $m =  2^{2^{16} + 1} $.
Since $2^{16} = A_11$, we may re-write $m$ as $A_0 (A_11 + 1)$.
Note that $A_11 < m < A^2_0 A_11 = A_12 < A_1A_11 = A_20$.
From this we obtain $A_10 < m <A_20$ and hence $a_1 = 1$, as well as $A_11<m<A_12$ which yields $b_1 = 1$, so that $m_1 = A_11$.
From $A_0A_11<m < A_12 = A_0^2 A_11 < A_0^2A_1(A_11-1) = A_1A_11$ we obtain $a_2 = 0$.
Moreover, $A_0(2^{16} + 1) = m < A_0(2^{16} + 2)$ yields $b _2=2^{16} + 1$.
Thus we have that $m \nfpar 2 A_0(2^{16} + 1)$. We may re-write this as $m \nfpar 2 A_0( A_{A_00} A_00 + A_00 )$; the reader may verify that all subexpressions are also in normal form.
\end{example}

\begin{remark}
Definition \ref{defNF} is not the only reasonable definition for normal forms, but it has advantages over other obvious candidates; see Appendix \ref{secConc}.
\end{remark}

Next we define the base-change operation based on our normal forms.

\begin{definition}\label{defBCH}
Given $k,\ell \geq 2$ and $m\in \mathbb N$ we define the base change operation $\bch m {\bases k \ell}$ recursively by setting $\bch 0 {\bases k \ell} = 0$ and, for $m \equiv _ k  A_{a} b + c$ setting
\[\bch m {\bases k \ell} =  A_{\bch {a} {\bases k \ell}} \bch b {\bases k \ell} +  \bch c{\bases k \ell}.\]
\end{definition}

Sometimes we abbreviate  $  \bch m{\bases {k}\ell}$ by $\bch m{\ell}$, in which case it is assumed that $k=\ell-1$ unless a different value for $k$ is specified.

\begin{example}
Let us write $A_xy$ for $A_x(2,y)$ and $B_xy$ for $A_x(3,y)$.
Recall from Example \ref{exHard} that $ 2^{16}+1 \gknf  A_11 + 1$.
Then, $\bch{  (2^{16}+1) }{\bases 23} =  B_11+1$.
We have that $B_11 = B_0^61 $, and hence
$\bch{  (2^{16}+1) }{\bases 23} = 3^{3^{3^{3^{27}}}} + 1$.
\end{example}

With this we are ready to define the Ackermannian Goodstein process.

\begin{definition}
Given a natural number $m$ we define a sequence $(\good mi)_{i< {\xi}}$ where ${\xi} \leq \omega$, by the following recursion.
\begin{enumerate}

\item $\good m0 = m$;

\item if $\good m{i} > 0$ then $\good m{i+1} = \bch {\good m i} { i +  3} - 1$;

\item if $\good mi = 0$ then ${\xi} = i+1$; if no such $i$ exists, $\xi = \omega$.

\end{enumerate}
The sequence $(\good mi)_{i<\xi}$ is the {\em Ackermannian Goodstein sequence starting on $m$.} 
\end{definition}

Our main results are then the following.

\begin{theorem}\label{theoMain}
Given any $m\in \mathbb N$ there is $i \in \mathbb N$ such that $\good mi = 0$.
\end{theorem}

\begin{theorem}\label{theoIndep}
Theorem \ref{theoMain} is not provable in $\atr$.
\end{theorem}

The rest of this article will be devoted to proving these results.
The proofs use properties of the ordinal $\Gamma_0$ reviewed in Appendix \ref{secGamma}.

\section{Properties of Normal Forms}

In this section we elaborate on the properties of normal forms as given by Definition \ref{defNF}.
Fix $k$ and write $A_xy$ for $A_x(k,y)$.
The intuition is that we obtain the normal form of $m$ by `sandwiching' it in smaller and smaller intervals, so that
\begin{equation}\label{eqIntervals}
  \big [ A_{a_1} (b_1), A_{a_1} (b_1+1)  \big ) \supsetneq \ldots \supsetneq \big [ A_{a_n} (b_n) , A_{a_n} (b_n+1)  \big ) \ni m.
\end{equation}
\longversion{
Recall that we have defined $m_0 = 0$ and $m_i = A_{a_i}b_i $ for $0<i\leq n$.
We have already remarked that $m_i < m_{i+1}$; in fact, the sandwiching indices and arguments are also strictly monotone.}
\shortversion{The following basic properties can readily be verified by the reader.}

\begin{lemma}\label{lemmSandwichingMonotone}
Let $(A_{a_i}b_i)_{i=1}^n$ be the sandwiching sequence for $m$ and $(m_i)_{i\leq n}$ be its sandwiching values.
\begin{enumerate}

\item\label{itSandwichingBM} If $0\leq i < n$ then $m_i \leq b_{i+1} < m_{i+1}$.

\item\label{itSandwichingA} If $0<i < n$ then $a_i>a_{i+1}$.

\end{enumerate}
\end{lemma}

\longversion{
\proof
For the first item note that $A_{a_{i+1}} m_i \leq m < A_{a_{i+1}} (b_{i+1} + 1)$, so, since $A_{a_{i+1}}$ is strictly increasing, $m_i \leq b_{i+1}$.
Moreover the inequality $A_xy >y$ yields $m_{i+1} = A_{a_{i+1}}b_{i+1}> b_{i+1}.$

For the second item, if we had that $a_i\leq a_{i+1}$ then we would also have that $A_{a_i}m_i \leq A_{a_{i+1}}m_i \leq m< A_{a_i} (b_i + 1)$, which is impossible since $m_i \geq b_i + 1$ by the first item.
\endproof
}

Thus we have that $A_{a_i} b_i < A_{a_{i+1}} b_{i+1}$, so in order to see that \eqref{eqIntervals} holds it remains to check that also $A_{a_{i}} (b_i+1) \geq A_{a_{i+1}} (b_{i+1}+1)$.

 \begin{lemma}\label{lemmSandwiching}
 Let $(A_{a_i}b_i) _{i = 1}^n$ be the $k$-sandwiching sequence for some $m >0$.
 Then, for any $i \in [1,n)$,
\begin{enumerate}
 

\item\label{itSandwichingTwo} $b_{i+1} < \min \{   A_{ a _{i+1}} ^{ k - 1} A_{{a}_{i+1} + 1} ( m _i - 1)  , A_{ a _{i} - 1} ^{ k - 1} m_i\}$;

\item\label{itSandwichingThree} $A_{  {a} _{i+1}} (b_{i+1} + 1)  \leq A_{{a} _{i+1} + 1} m_i $, and

\item\label{itSandwichingFour} if $i>0$ then
$A_{{a}_{i+1}} (b_{i+1} + 1)  \leq A_{{a}_{i}} (b_{i} + 1) $.
\end{enumerate}

 \end{lemma}
 
 \proof   \eqref{itSandwichingTwo}
That $b_{i+1} < A_{ a _{i+1}} ^{ k - 1} A_{{a}_{i+1} + 1} ( m _i - 1)$ follows from
$A_{ a _{i+1}} ^{ k } A_{{a}_{i+1} + 1} ( m _i - 1) = A_{a_{i+1}+1} m_i > m \geq A_{a_{i+1}} b_{i+1} $
and the monotonicity of $A_{a_{i+1}}$.
 
For the other inequality, if $a_{i+1} < a_i-1$ then the definition of $a_{i+1}$ yields $m < A_{a_i-1}m_i$ and there is nothing to prove. Otherwise, by Lemma \ref{lemmSandwichingMonotone}.\ref{itSandwichingA} we would have that $a_{i+1} = a_i-1 $.
But then we observe that
$A_{a_{i+1}}^k A_{a_i}b_i =  A_{a_i}(b_i+1) > m \geq A_{a_{i+1}} b_{i+1}$,
so that $ b_{i+1} <  A_{a_{i+1}}^{k-1} A_{a_i}b_i = A_{a_{i+1}}^{k-1}m_i$.
\medskip

\noindent \eqref{itSandwichingThree}
By item \eqref{itSandwichingTwo},
$A_{a _ {i+1} + 1 } m_i = A_{  {a} _ {i+1} }  A_{  {a} _ {i+1} }^{k-1} A_{a_ {i+1} + 1 } ( m_i - 1)   \geq A_{  {a}_ {i+1} } (b _ {i+1} + 1)$.

\noindent \eqref{itSandwichingFour}
Once again by item \eqref{itSandwichingTwo},
$A_{a_i} (b_i+1) 
 = A_{a_i - 1}^k m_i
  \geq A_{a_{i+1}} (b_{i+1} + 1)$.
 \endproof

The following two propositions provide the basic techniques for computing the normal form of $w$ assuming $m\gknf A_ab+c$, where $w$ is ``not too different'' from $m$.

\begin{proposition}\label{propNewNormF}
Suppose that $m $ has $k$-sandwiching sequence $(A_{a_i}b_i)_{i=1}^n$ and $w$ and $1\leq j\leq n$ are such that
\[A_{a_j}b_j \leq w < A_{a_j} (b_j+1).\]
Then, $(A_{a_i}b_i)_{i=1}^j$ is an initial segment of the $k$-sandwiching sequence for $w$.

If moreover $A_0m_{j} > w$, then $(A_{a_i}b_i)_{i=1}^j$ is the full $k$-sandwiching sequence for $w$ and hence $w\gknf A_{a_j}b_j + c$ for suitable $c$.
\end{proposition}

\proof
Assume that $A_{a_j}b_j \leq w < A_{a_j} (b_j+1)$ and let $(A_{d_i}e_i)_{i=1}^r$ be the $k$-sandwiching sequence for $w$.
In view of Lemmas \ref{lemmSandwichingMonotone} and \ref{lemmSandwiching} we have for $i<j$ that
$A_{a_{i+1}} m_i  \leq A_{a_{i+1}} b_{i+1} \leq A_{a_j}b_j
\leq w < A_{a_j} (b_j+1) \leq A_{a_{i+1}} (b_{i+1}+1) \leq A_{a_{i+1} + 1} m_i$,
so that $(A_{a_i}b_i)_{i=1}^j$ satisfies the recursion of Definition \ref{defNF} applied to $w$.
It follows that $r\geq j$ and for $0<i\leq j$ that $a_i = d_i$ and $b_i = e_i$.
If moreover $A_ 0 m_j > w $ then the sandwiching halts and $r=j$, witnessing that $w\gknf A_{a_j}b_j + c$ for some $c$.
\endproof

As a special case, if $1\leq j$ we can set $w= m_j$ and see that $(A_{a_i}b_i)_{i=1}^j$ is the sandwiching sequence for $m_j$, hence $m_j \gknf A_{a_j}b_j$.
Note that $(A_{a_i}b_i)_{i=1}^0$ is also the sandwiching sequence for $m_0$, since $0 = m_0$ has empty sandwiching sequence.

\longversion{
Proposition \ref{propNewNormF} identifies numbers that share a portion of their sandwiching sequence.
Next, we characterize some cases where new values are added to the sandwiching procedure.
}

\begin{proposition} \label{propNormForm}
Let $m \gknf A_a b$ with sandwiching sequence $(A_{a_i}b_i)_{i=1}^n$, $j \leq n$ and $w =A_de$. Suppose that
\begin{enumerate}[wide]

\item $A_dm_j \leq w <  A_{d+1}m_j $, and

\item if $j> 0$ then $ w < A_{a_j} ( b_j + 1)$.

\end{enumerate}
\noindent Then, $w \gknf A_de$ and has sandwiching sequence $(A_{d_i}e_i)_{i=1}^{j+1}$, where $d_i=a_i$ and $e_i=b_i$ for $i\in[1,j]$, $d_{j+1} = d $, and $e_{j+1} = e$.
\end{proposition}

\proof
Let $( A_{d _i } e_i )_{i=1}^r$ be
the sandwiching sequence for $ w $.
%
%
%
We claim that $r\geq j$ and for $1\leq i \leq j$, $d_ i = a_i$ and $e_i = b_i$.
When $j=0$ this is vacuously true, otherwise this follows from the inequality
$A_{a_j}b_j = m_j < A_dm_j \leq w < A_{a_j}(b_j+1) $
and Proposition \ref{propNewNormF}.
In particular, $w_{j} = m_j$.
Now, the inequality
$ A_{d} m_j \leq w < A_{d+1} m_j$
yields $ d_{j+1} = d $, and since $w=A_de<A_0A_de$ the sandwiching halts and $w \equiv_k A_{d} e $.
\endproof

\begin{corollary}\label{corSmallIndex}
Let $m \gknf A_a b$ with sandwiching sequence $(A_{a_i}b_i)_{i=1}^n$, $1\leq x < k $, and $j\leq n$ be such that if $j \geq 1$ then $d<a_j$.
Then, $A_d^x m_j$ is in normal form.
\end{corollary}

\longversion{
\proof
It is easy to check that the required inequalities for Proposition \ref{propNormForm} hold.
\endproof
}

We omit the proof, which consists of checking that the required inequalities for Proposition \ref{propNormForm} hold.
Now let us introduce some additional notation for sandwiching sequences.
If $(x_i)_{i=n_0}^n$ is any sequence with $0\leq n_0<n$, we denote its second-to-last element by $x_{-1} $, i.e.~$x_{-1}:=x_{n-1}$.
If $m>0$ and $(m_i)_{i=0}^n$ are its sandwiching values, we define $\penum m = \penu{m}$; note in particular that if $(A_{a_i} b_i)_{i=1}^n$ is the sandwiching sequence for $m$ and $n>1$ then $\penum m =A_{\penu a}\penu b $.

\begin{lemma}\label{lemmNormForm}
Let $ a,b,  e$ be natural numbers and let $m \gknf A_{a} b$ and $w =  A_{a-1} e $.
Then, $w\gknf A_{a-1} e $ if either
\begin{enumerate}

\item\label{itSmall} $ A_{a-1} \penum{m} \leq w  <  A_{a} \penum{m}$,
or

\item\label{itBig} $A_{a-1} m \leq w < A_{a} (b+1) $.
\end{enumerate}

\end{lemma}

\proof
Let $(A_{a_i}b_i)_{i=1}^n$ be the sandwiching sequence for $m$.
In the first case, if $n=1$ this is an instance of Proposition \ref{propNormForm}. Otherwise, note that $\penum m=A_{a_{-1}}b_{-1} $, and since $ a<a_{-1}$ we obtain
$A_{a-1} \penum m \leq w < A_{a} \penum m <A_{a_{-1} - 1}^k A_{a_{-1}} b_{-1} = A_{a_{-1}} (b_{-1}+1)$. 
In the second we instead have
$A_{a-1} m \leq w < A_{a } (b +1) < A_a A_ab = A_{a } m$.  
Thus in both cases we can apply Proposition \ref{propNormForm}.
\endproof

Normal forms can be divided into the cases where $b=\penum m$ and $b>\penum m$, as each behaves in a different manner.
The following lemma makes this precise.

\begin{lemma}\label{lemmNormBCases}
Let $( A_{{a}_i}b_i)_{i=1}^n$ be the $k$-sandwiching sequence for $m \equiv_k A_{a} b $.
Then, exactly one of the following cases occurs:
\begin{enumerate}[label=(\alph*)]

\item\label{itBCasesOne} $n=1$ and $ \penum{m} = b = 0$,

\item\label{itBCasesTwo} $n>1$ and $0< \penum{m} = b \equiv_k   A_{\penu{a}} \penu{b} $, or

\item\label{itCaseThree} $n\geq 1$ and $  \penum{m} < b \equiv_k   A_d e + s $ with $d \leq \penu{a}$ and either $d \leq a$ or $s>0$.

\end{enumerate}
\end{lemma}

\proof
Assume that \ref{itBCasesOne} fails, so that $ n > 1$ or $b>0$.
Let $( A_{d_i} e_i)_{ i = 1}^{r}$ be the $k$-sandwiching sequence for $b$.
Note that $( A_{{a}_i}b_i)_{ i = 1}^{n-1}$ is a (possibly empty) initial segment of this sequence: if $n=1$ this holds trivially, while for $n>1$ we have that $A_{\penu{a}} \penu{b} \leq b \leq A_{\penu{a}} ( \penu{b} + 1)$, hence we can apply Proposition \ref{propNewNormF}.
Thus $r \geq n-1$ and for $1\leq i < n$ we have that $d_i=a_i$ and $e_i = b_i$. Now, consider two cases.
\begin{enumerate}[label*={\sc Case \arabic*},wide, labelwidth=!, labelindent=0pt]

\item ($ r =  n-1$). This would witness that $ b\equiv_k A_{{a}_{ -1}}\penu{b}+s $ for some $s$. 
Hence ${d} = {a}_{ -1}   $ and $\penu{b}= e$, so that either $ b = \penum{m}$ and \ref{itBCasesTwo} holds, or $b>\penum{m}$ and $s>0$, so that \ref{itCaseThree} holds.
\medskip

\item ($r\geq n$). In this case we claim that $  {d} _n \leq {a}$. Indeed, $A_{a+1} A_{\penu d} \penu e = A_{a+1} \penum m > A_{a} b > b$.
Since $A_{{d}_n} A_{\penu a}	 \penu e \leq b$, we conclude that $  {a} \geq d_n \geq d_r = d$. This shows that $ b \equiv_k   A_{d_r} e_r + s $ for some $s$, and \ref{itCaseThree} holds.\qedhere
\end{enumerate}
\endproof

\begin{lemma}\label{lemmNormFormMinusOne}
Let $a,b,c,m,w $ be natural numbers such that
\[ A_ab + c =  w < \min \{  A_a(b+1), A_0 A_ab\}. \]
Then, $w \gknf A_ab + c$ if either
\begin{enumerate}

\item\label{itBMinusOne} $m \gknf A_{a} (b+1) > A_a \penum{m}$, or

\item\label{itCMinusOne} $m \gknf A_{a}  b + c'  $ for some $c'$.

\end{enumerate}
\end{lemma}

\proof
In the first case, it is not hard to check using Proposition \ref{propNewNormF} that
the sandwiching sequence for $ w $ is the same as that for $m$, except that the last item is replaced by $A_ab$.
In the second, $w$ has the same sandwiching sequence as $m$.
\endproof

Note that if $m \gknf A_a b + c$ we may still have that $c\geq A_ab$.
For such cases we define the {\em extended $k$-normal form} of $m$ to be $A_ab \cdot p +q$, in symbols $m\gnf  A_ab \cdot p +q$, if $m\gknf A_ab + c $ for some $c$, $m =   A_ab \cdot  p +q$, and $0 \leq q < A_a b$; note that $p$ and $q$ are uniquely defined.\footnote{The operations $A_xy$ and $\bch x{\bases yz}$ are always assumed to be performed before multiplication.}
If $m\gnf  A_ab \cdot p +q$ and $d = A_ab$ we write $m \cnf d\cdot p + q$ and call it the {\em simplified $k$-normal form} of $m$.
\shortversion{We state two useful corollaries of Lemma \ref{lemmNormFormMinusOne}. Both are proven by checking \eqref{itCMinusOne}; in the first we set $c=0$ and $c' = A_ab\cdot(p-1)$, and in the second $c =  d \cdot p + q$ and $c' = d\cdot (p-1)   + q$.
In both cases, the required inequality $w <\min \{A_a(b+1), A_0A_ab \}$ is not hard to verify.
}

\begin{corollary}\label{corMultK}
If $A_ab$ is in normal form and $0<p<k$ then $A_ab \cdot p$ is in extended normal form.
\end{corollary}

\longversion{
\proof
This follows from Lemma \ref{lemmNormFormMinusOne}.\ref{itCMinusOne} and the inequality
$A_ab \leq pA_ab <\min \{A_a(b+1), A_0A_ab \}$.
To see this, note that $pA_ab < A_a(b+1)$ follows from $A_a(b+1) = kA_ab$ if $a=0$ and $A_a(b+1) \geq A_0^2 A_ab \geq A_0(A_ab+1) = k A_ab$ if $a>0$.
Meanwhile,
$pA_ab < k^{A_ab} = A_0A_ab  $ follows from the inequality
$(x+1)^y > xy 
$
that holds for all natural numbers $x,y$.
\endproof
}

\longversion{
Note that in Lemma \ref{lemmNormFormMinusOne}.\ref{itCMinusOne} it is sufficient (but not necessary) to have $w < m$, since in this case we already have $m < \min \{A_0 m,A_a(b+1)\}$.
With this in mind, we obtain the following corollary of Lemma \ref{lemmNormFormMinusOne}
}

\begin{corollary}\label{corNFC}
If $m \cnf  d \cdot (p+1)  + q$ with $p>0$ then $  d \cdot p +q$ is in simplified normal form.
\end{corollary}

\longversion{
\proof
Immediate from Lemma \ref{lemmNormFormMinusOne}.\ref{itCMinusOne} by setting $c =  d \cdot p + q$ and $c' = d\cdot (p-1)   + q$.}
\endproof

As we will see in the next section, if $m \gknf A_a b$ it is usually not easy to compute the normal form of $m-1$.
However, the case where $a=0$ is somewhat simpler.

\begin{lemma}\label{lemmZeroNF}
Fix $ k \geq 0 $ and write $A_xy  = A_x (k,y)$.
Suppose that $m\gknf A_0b$.
\begin{enumerate}

\item \label{itZeroBBig} If $b > \penum{m}$ then $m-1 \cnf  A_0 (b-1) \cdot(k-1) + (k^{b-1} - 1)$.

\item \label{itZeroBSmall} If $b = \penum{m} > 0 $ then $m-1 \cnf  b \cdot p + q$ for some $p < k^{b - 1}$ and some $q< b$.

\end{enumerate}
\end{lemma}

\proof
We prove \eqref{itZeroBSmall}; item \eqref{itZeroBBig} is similar.
Suppose that $m\gknf A_0  b$ and let $(A_{a_i}b_i)_{i=1}^n$ be the sandwiching sequence for $m$.
By Lemma \ref{lemmNormBCases}, $n>1$ and $  b = A_{\penu{a}} \penu{b} $.
It is then easy to check using Proposition \ref{propNewNormF} that $ m-1 \gknf A_{\penu a}\penu b +c$ for some $c$, hence $m-1 \cnf  b \cdot p + q$ for some $p$, $q$.
Since $\penu{a} > a_n $ we must have that $\penu a\geq 1$.
It follows that $ b \geq A_1 0 > k$, so that if $p\geq k^{b - 1}$ we would have that $bp > k^{b - 1} \cdot k = A_0  b$, which is impossible as $bp \leq m-1$.
That $q < b$ follows from the definition of simplified normal forms.
\endproof

\section{Left and Right Expansions}

\longversion{
In the previous section we have shown that, under certain conditions, if $m \gknf A_ab + c$ then for suitable $b' \leq b$ and $c'\leq c$, $A_a b' + c'$ is in normal form.
However, there are also cases where such expressions are not in normal form: for example, we may only conclude that $A_a(b-1)$ is only in normal form if $b>\penum m$.
In this section we will show how such normal forms can be found by expanding $A_a$ according to its recursive definition.
First we expand it on the right in order to compute the normal form of $A_a(b-1)$ when $b = \penum m$.
Throughout this section we once again fix $ k > 1 $ and write $A_xy  $ for $A_x (k,y)$.
Below, when $s=b+1$ recall that by definition $A_a(-1)=1$.
}

\shortversion{
In this section we develop operations for computing normal forms using the information that $m \gknf A_ab + c$, but where the new normal form may vary drastically from that of $m$.
We begin with {\em right expansions} of $A_ab$.
}

\begin{lemma}\label{lemmNormFormZero}
Let $m = A_{a} b $ with ${a}, b>0$, and
let $s \in [ 1 , b + 1]$.
Then:

\begin{enumerate}
\item\label{itNFZPenum} $A_{a} b  = A^{sk}_{ a - 1}A_{a} (b -s ).$

\item\label{itNFZLast} Let $\ell \in [1,k]$ and $ c = A^\ell _{a-1} A_{a} (b -s  )$. If $m \equiv_k A_{a} b $, $b = \penum{m}$, and
\[ A_{ a - 1} b \leq c < A_ab,\]
it follows that $c $ is in normal form as written.

\end{enumerate}
\end{lemma}

\proof
The first claim follows by a simple induction on $s$ and the definition of $A_{a}$ and the second is an instance of Lemma \ref{lemmNormForm}.\ref{itSmall}.
\endproof

The following is a variant of the ordinal predecessor function of Cichon \cite{Cichon}, which works by expanding $A_ab$ on the {\em left,} useful for computing the normal form of $m-1$.

\begin{definition}\label{defLeftExp}
Let $A_{a} b$ be in normal form with $a>0$ and define a sequence $c_0,\ldots, c_{a}$ by recursion as follows:
\begin{enumerate}

\item $ c_0 =
A_{a} (b - 1)$;

\item $c_{i} =
A_{ a - i } ( A^{ k - 1} _{ a - i } c_{i-1} - 1)
$ if $i>0$.

\end{enumerate}
\noindent We call the sequence $(c_i)_{i\leq a}$ the {\em left expansion sequence} for $A_{a} b$.

\end{definition}

The following chain of inequalities summarizes some basic properties of left expansion sequences and can be easily checked using induction on $i$.

\begin{lemma}\label{lemmNFLeftExpBasic}
If $m \equiv _k A_{a} b$ with $a>0$, $(c_i)_{i\leq a}$ is the left expansion sequence for $m$, and $ i \in [1,a]$, then
\[b < c_{i-1} <c_{i} < m = A_{a-i}^{k} c_{i-1} < A_{a-i+1} c_{i-1}.\]
\end{lemma}

In fact, the inequality $c_i<m$ can immediately be sharpened.

\begin{corollary}\label{corKCI}
Suppose that $m\gknf A_ab$ with $a>0$ and let $(c_i)_{i\leq a}$ be the left expansion sequence for $m$.
Then, if $i\leq a$, it follows that $m \geq k c_i $, with equality holding if and only if $i=a$.
\end{corollary}

\proof
For $i=a$, note that
$m= A_{0}^{k} c_{a-1} = k A_0(A_0^{k-1} c_{a-1} -1) = k c_a$. 
Otherwise, $kc_i < k c_a = m$.
\endproof

\longversion{
Next we show that the elements of a left expansion sequence and their displayed subterms are in normal form.
}

\begin{lemma}\label{lemmNFLeftExp}
Let $m \equiv _k A_{a} b$ with $a>0$.
\begin{enumerate}

\item If $ 0 < \ell < k$ and $0 <  i \leq a$ then $A_{a-i}^\ell c_{i-1} $ is in normal form.

\item If $ i \leq a$ and either $i > 0$ or $b>\penum{m}$, then $c_i$ has normal form as written in Definition \ref{defLeftExp}.

\end{enumerate}
\end{lemma}

\proof
Proceed by induction on $i$, considering several cases.
\begin{enumerate}[label*={\sc Case \arabic*},wide, labelwidth=!, labelindent=0pt]

\item ($i=0$ and $b>\penum{m} $).
The first claim does not apply with $i=0$ and the second is an instance of Lemma \ref{lemmNormFormMinusOne}.\ref{itBMinusOne}, using the condition $b>\penum m$.
\medskip

\item ($i>0$). We will prove both claims uniformly. In order to do this, let $d$ be either $A_{a-i}^{\ell-1} c_{i-1}$ for the first claim or $A_{a-i}^{k-1} c_{i-1} - 1$ for the second.
Note that it follows from Lemma \ref{lemmNFLeftExpBasic} that $A_{a-i} d < m $.
We will consider two sub-cases.
\medskip

\begin{enumerate}[label*= .\arabic*,wide, labelwidth=!, labelindent=0pt]

\item ($b=\penum{m}$ and $i = 1$).
This is an instance of Lemma \ref{lemmNormFormZero}.\ref{itNFZLast}.
\medskip

\item ($b>\penum{m}$ or $i>1$).
By the induction hypothesis we have that $c_{i-1} \gknf A_{a-i+1} e$ for some $e$, where $e=b-1$ if $i=1$ and $e=A^{k-1}_{a-i+1} c_{i-2} -1$ if $ i> 1 $.
In either case $A_{a-i+1}(e+1) = m$, in the first case since $m=A_{a }b$ and in the second by Lemma \ref{lemmNFLeftExpBasic}.
Moreover, $c_{i-1}\leq d$, so
$A_{a-i }  c_{i-1} \leq A_{a-i}d <  m = A_{a-i+1}(e+1)$. 
Thus we may use Lemma \ref{lemmNormForm}.\ref{itBig} to conclude that $A_{a-i}d$ is in normal form.\qedhere
\end{enumerate}
\end{enumerate}
\endproof

With this we can describe the normal form of $m-1$ when $m\gknf A_ab$.

\begin{lemma}\label{lemmNormFormMinus}
If $m\equiv_k A_a b$ with $a>0$ and left expansion sequence $(c_i)_{i\leq a}$ then
\[m-1 \cnf (k-1) \cdot c_a +  (c_a - 1).\]
\end{lemma}

\proof
Let $d = A_0^{k-1} c_{a-1}-1$.
We know that $c_a \gknf  A_0 d$ and
$A_0 d  = c_a < m-1 < m = A_0 ( d+1 ) \leq A_0 A_0 d$. 
Thus Lemma \ref{lemmNormFormMinusOne}.\ref{itCMinusOne} implies that $ m-1 \gknf A_0 d+ s $ for some $s$, hence $m-1\cnf  c_a \cdot p + q$ for some $p,q$ with $q<c_a$.
But by Corollary \ref{corKCI},
$m-1 = kc_a - 1 =  c_a \cdot (k-1)  + (c_a - 1) $, 
so that $p=k-1$ and $q = c_a - 1$.
\endproof

\section{The Base Change Operation}

Next we elaborate on the base change operation as given by Definition \ref{defBCH}.
We begin by extending it to base $\omega$ as follows:
\begin{enumerate}
\item If $m=0$ then $\pk m =0$.
\item If $m\gknf A_a b + c$ then 
$\pk m = \ot_{\pk a}\pk b+ \pk c$.

\end{enumerate}
Here, $\ot_\alpha$ denotes the fixed point-free Veblen functions, reviewed in Appendix \ref{secGamma}.
We regard $\ot_\alpha$ as an analogue of the Ackermann function with base $\omega$.
To stress this, we set $A_\alpha(\omega,\beta) : = \ot_\alpha \beta$, so that Definition \ref{defBCH} uniformly extends to $\ell \leq \omega$.
The following is shown by induction on $p$ using Corollary \ref{corNFC}.

\begin{lemma}\label{lemmBCHcoeff}
If $m \gnf A_a(k,b) \cdot p + q $ and $1 < k\leq \lambda \leq \omega $ then
\[\bch m {\bases k\lambda} =  A_{\bch a {\bases k \lambda}} (\lambda, \bch b {\bases k \lambda}) \cdot p + \bch q {\bases k\lambda}.\]
\end{lemma}

Next we prove that the base change operation is strictly increasing, which will be a crucial ingredient in the proofs of our main results.

\begin{proposition}\label{propMonotonicity}
If $n < m$ and $1 < k\leq \lambda \leq \omega $ then $\bch n{\bases k \lambda } < \bch  m {\bases k \lambda } $.
\end{proposition}

\proof
By induction on $m$.
In this proof we write $\bch x\lambda $ instead of $\bch x{\bases k\lambda}$.
We remark that if $x<m$ our induction hypothesis yields $x\leq \bch x\lambda$, since the map $ x\mapsto \bch x  \lambda$ is strictly increasing below $m$.
Without loss of generality we may assume that $n = m - 1$.
Let $A_xy = A_x (k,y)$ and $B_\xi\zeta = A_\xi (\lambda, \zeta )$. Write $m \equiv_k  A_{a} b + c$, so that $\bch m \lambda =  B_{\bch {a} {\lambda}} \bch b\lambda   + \bch c \lambda$.
\shortversion{We then consider several cases.
}

\begin{enumerate}[label*={\sc Case \arabic*},wide, labelwidth=!, labelindent=0pt]

\item ($a=b=c=0$). Then $m=1$ so $\bch n\lambda = n = 0   <\bch m\lambda$.
\medskip

\item ($c>0$).
In this case we have by Lemma \ref{lemmNormFormMinusOne}.\ref{itCMinusOne} that $m-1 \equiv_k A_{a} b + (c-1)$, so that
$\bch {(m-1)} \lambda = B_{\bch a \lambda } \bch b \lambda + \bch {(c - 1)} \lambda  <^{\text{\sc ih}} B_{\bch a\lambda} \bch b \lambda + \bch { c } \lambda = \bch m\lambda$.
\medskip

\item ($a=c=0$ and $\penum{m} < b$).
Write $b=d+1$. Lemma \ref{lemmZeroNF}.\ref{itZeroBBig} then yields $m-1\gnf  A_0d \cdot (k-1) + (k^d - 1)$, hence
\begin{align*}
\bch m\lambda & = B_{\bch 0\lambda} \bch{b }\lambda
 \geq^{\text{\sc ih}}  B_{0} (\bch {d} \lambda +1 )
> B_{0}  \bch {d} \lambda   \cdot k
 = B_{0}  \bch {d} \lambda   \cdot (k-1) + B_{0}  \bch {d} \lambda \\
&
  >^{\text{\sc ih}} B_{0}  \bch {d} \lambda  \cdot (k-1)+\bch {(k^d-1 )}\lambda = \bch{(m-1)}\lambda.
\end{align*}

\item ($a=c=0$ and $\penum{m} = b$).
Use Lemma \ref{lemmZeroNF}.\ref{itZeroBSmall} to see that $m-1 \cnf  b \cdot p + q$ for some $p < k^{b - 1}$ and some $q< b$.
Then,
\begin{align*}
\bch m\lambda  & = B_0 \bch b\lambda 
 \geq  \bch b \lambda \cdot k^{\bch b\lambda  -1} + \bch b \lambda\\
& \geq^{\text{\sc ih}}  \bch b \lambda \cdot k^{b-1} + \bch q \lambda
 \geq   \bch b\lambda \cdot p + \bch  q \lambda
=\bch {(m-1)} \lambda,
\end{align*}
where for $\lambda<\omega$ the first inequality follows from the fact that $x^{y-1} y + y \leq (x+1)^y
$ whenever $x,y$ are natural numbers with $x>0$.
\medskip

\item ($a>0$ and $c= 0$).
Let $c_0,\ldots, c_{a}$ be the left expansion sequence for $m$.
If $\lambda < \omega$ let $ \delta_0,\ldots, \delta_{\bch a \lambda }$ be the left expansion sequence for $\bch m \lambda$, and if $\lambda = \omega$ define $\delta_i = \bch m\lambda$ for all $i$.
If $\lambda<\omega$ note by the induction hypothesis that $a-i \leq \bch {(a-i)}\lambda\leq \bch a \lambda - i$ whenever $i\leq a$.
We claim that for each $i \leq a$, $\bch {c_i} \lambda < \max \{2,  \delta_i\}$.
The theorem will then follow, since by Lemma \ref{lemmNormFormMinus} we see that
$m-1 \cnf  c_a \cdot (k-1)  +  (c_a - 1),$
hence
$ \bch{(m-1)} \lambda
=   \bch{ c_a }\lambda \cdot (k-1) + \bch{ (c_a - 1) } \lambda  <^{\text{\sc ih}} \bch{ c_a }\lambda \cdot  k < \bch m \ell$,
where the last inequality follows from Corollary \ref{corKCI} if $\lambda < \omega$.

We proceed to prove that $\bch {c_i} \lambda < \max \{2,  \delta_i\}$ by a secondary induction on $i$. Let us first do the inductive step:
$\bch{c_{i+1}}\lambda= \bch{ A_{a-i-1} ( A^{ k - 1} _{( a - i - 1 ) } c_i - 1) }\lambda 
 =B_{\bch{( a - i - 1 ) } \lambda } \bch{ ( A^{ k - 1} _{( a - i - 1 ) } c_i - 1) } \lambda 
 <  B_{\bch{( a - i - 1 ) } \lambda}    \bch{ A^{ k - 1} _{( a - i - 1 ) } c_i}\lambda  
 =    B^{ k } _{\bch{( a - i - 1 ) } \lambda} \bch {c_i}\lambda 
< \delta_{i+1} $,
where the first inequality follows from $A^{ k - 1} _{( a - i - 1 ) } c_i < m$ and the main induction hypothesis and the last from $ B^{ k } _{\bch{( a - i - 1 ) } \lambda} \bch {c_i}\lambda<  B _{\bch a\lambda - i -1}   \big ( B^{ k } _{\bch a\lambda - i -1} \delta_i - 1\big ) = \delta_{i+1} $ when $\lambda < \omega$ and from $\bch m\ell$ being closed under $B^{ k } _{\bch{( a - i - 1 ) }\lambda}$ when $\lambda = \omega$.

Thus it remains to prove that $\bch {c_0} \lambda \leq \delta_0$.
\medskip

\begin{enumerate}[label*= .\arabic*,wide, labelwidth=!, labelindent=0pt]
\item ($b=0$). We have that $ \bch {c_0}\lambda = c_0  = 1 < 2$.
\medskip

\item ($b> \penum{m}$).
We have that
$\bch {c_0} \lambda = B_{\bch {a}\lambda} \bch { (b-1)} \lambda < \delta_0 $,
where if $\lambda<\omega$ we use the induction hypothesis to see that $ \bch { (b-1)} \lambda \leq \bch { b } \lambda - 1 $ and if $\lambda = \omega$ we use $\bch { (b-1)} \lambda < \bch b\lambda $ and the monotonicity of $B_{\bch a \lambda}$.
\medskip

\item\label{itCritPropMon} ($ b = \penum{m}>0$).
We use Lemma \ref{lemmNormFormZero} to see that for $ s \in [1,b]$,
$A_{a} (b-1) = A^{k(s-1)}_{{a} - 1}  A_{a}(b-s) $.
Since $A_{a} (-1)  = 1 \leq b$ we have that there is a least $t \in [2 , b + 1] $ such that
$A_{a}(b-t ) \leq b $.
Similarly, there is a greatest $r < k$ such that
$u  := A^r_{{a} - 1} A_{a}(b-t ) \leq b  $.
Note that $A_{a-1} u>b$, and hence $A^2_{a-1} u > A_{a-1}b = A_{a-1}\penum m$.
Hence by Lemma \ref{lemmNormFormZero}.\ref{itNFZLast} we have that
$A_{{a} - 1}^{v}  u  $
is in normal form whenever
$1 < v \leq k (t-2) + k - r .$
Similarly, by Corollary \ref{corSmallIndex}, $ A_{{a} - 1} b$ is in normal form, since by assumption $b=\penum{m}$.
Then,
\begin{align}
\nonumber
\bch {c_0} \lambda & = \bch { A^{k(t-2) + k-r-1}_{a-1} A_{{a} -1 } u  } \lambda
=  B^{k(t-2) + k-r-1}_{\bch {(a-1)} \lambda } \bch {   A_{{a} -1 } u  } \lambda\\
&\nonumber \leq^{\text{\sc ih}} B^{k(t-2) + k -1}_{\bch {(a-1)} \lambda } \bch {   A_{{a} -1 } b   } \lambda
\nonumber = B^{k(t-2) + k - 1  }_{\bch {(a-1)} \lambda }  B_{\bch {(a-1)} \lambda } \bch { b } \lambda \nonumber = B^{k(t-1)  }_{\bch {(a-1)} \lambda } \bch { b } \lambda.
\end{align}
If $\lambda = \omega$ then $\bch {(a-1)} \lambda < \bch a \lambda $ yields $B^{k(t-1)  }_{\bch {(a-1)} \lambda } \bch { b } \lambda < B _{\bch {a} \lambda } \bch { b } \lambda = \bch m \lambda$.
If $\lambda < \omega$, then
$
B^{k(t-1)  }_{\bch {(a-1)} \lambda } \bch { b } \lambda \nonumber < B^{k(t-1)  }_{\bch {(a-1)} \lambda } B^{t-1 }_{\bch {(a-1)} \lambda } B_{\bch a \lambda} (\bch { b } \lambda - t)
\nonumber \leq B^{\lambda(t-1)  }_{\bch {(a-1)} \lambda }  B_{\bch a \lambda} (\bch { b } \lambda - t)
\nonumber \leq B^{\lambda(t-1)  }_{\bch {a} \lambda - 1 }  B_{\bch a \lambda} (\bch { b } \lambda - t)
\nonumber = B_{\bch a \lambda} (\bch { b } \lambda - 1) = \delta_0$,
as claimed.\qedhere
\end{enumerate}

\end{enumerate}
\endproof

 \section{Normal Form Preservation}

Our goal now is to show that the base-change operation preserves normal forms: if $A_a(k,b)+c$ is in $k$-normal form and $\lambda \geq k$, then $A_{\bch a{\bases k\lambda}} (\lambda, \bch b{\bases k\lambda}) + \bch c{\bases k\lambda}$ is in $\lambda$-normal form.
Throughout this section we fix $2 \leq k < \lambda \leq \omega$, and write $A_xy = A_x(k,y)$, $B_\xi\zeta =A_\xi(\lambda,\zeta)$.

 \begin{proposition}\label{propNFP}
 If $m\gknf A_a b + c $ then $\bch m \lambda \equiv_{\lambda} B_{\bch a\lambda } \bch b\lambda + \bch c \lambda$.
 \end{proposition}

 \proof
Assume that  $m\gknf A_a b + c$.
 The proposition is immediate from Proposition \ref{propMonotonicity} when $\lambda = \omega$, since the only condition for $\bch m \omega = \ot_{\bch a {\bases k\omega}} \bch b {\bases k\omega} + \bch c{\bases k\omega}$ to be in normal form is for $\bch a{\bases k\omega}, \bch b {\bases k\omega}, \bch c {\bases k\omega} < \bch m {\bases k\omega} $, which holds since $a,b,c<m$.
 Thus we assume that $\lambda <\omega$ and without loss of generality that $\lambda = k+1$. Let $(A_{a_i} b_i)_{i=1}^n$ be the $k$-sandwiching sequence for $m$.
We will prove by induction on $m$ that $(B_{\bch {a_i} \lambda } \bch {b_i} \lambda )_{i=1}^n$ is the $\lambda$-sandwiching sequence for $\bch m \lambda$.

\begin{enumerate}[label*={\sc Case \arabic*},wide, labelwidth=!, labelindent=0pt]

\item ($c>0$).
In this case $m - 1 \equiv_k A_0 b + (c-1)$, so that by the induction hypothesis $(B_{\bch {a_i} \lambda } \bch {b_i} \lambda )_{i=1}^n$ is the $\lambda$-sandwiching sequence for $\bch {( m - 1)} \lambda$.
Consider two sub-cases.
\medskip

\begin{enumerate}[label*= .\arabic*,wide, labelwidth=!, labelindent=0pt]
\item ($a=0$). Since $B_0 B_0\bch{b }\lambda \geq B_0(\bch{b }\lambda + 1)$, to apply Proposition \ref{propNewNormF} we need only check that $B_0(\bch{b }\lambda + 1) > \bch m\lambda$.
Write
$m \gnf  A_0 b \cdot p + q = k^bp + q $
with $q<k^b$.
From $ m < A_0 (b +1) $ we obtain $p < k$.
However,
$B_0(\bch{b }\lambda + 1)
> B_0 \bch b \lambda \cdot k  
 \geq \bch {A_0 b} \lambda \cdot p + \bch {A_0 b} \lambda
 >  \bch {A_0 b} \lambda \cdot p  + \bch q \lambda 
  = \bch m\lambda$.
  \medskip


\item ($a>0$).
Since $ B_{\bch a\lambda} (\bch b\lambda + 1)  = B_{\bch a\lambda-1}^\lambda B_{\bch a\lambda} \bch b\lambda > B_0 B_{\bch a\lambda} \bch b\lambda $,
now it suffices to show $B_0 B_{\bch {a } \lambda } \bch {b } \lambda > \bch m \lambda$ to conclude that $(B_{\bch {a_i} \lambda } \bch {b_i} \lambda )_{i=1}^n$ is also the sandwiching sequence for $\bch m\lambda$.
We know that
$ A_0 A_{a } b  > m.$
Since $a  > 0$ we have that $A_0 A_{a } b $ is in normal form by Corollary \ref{corSmallIndex}, hence by Proposition \ref{propMonotonicity},
$B_0B_{\bch{a }\lambda}\bch {b } \lambda = \bch{A_0 A_{a } b } \lambda > \bch m \lambda$, 
as claimed.
\end{enumerate}
\medskip

\item ($c=0$).
This is the critical case in the proof.
By Proposition \ref{propNewNormF}, $\penum m $ has $k$-sandwiching sequence $(A_{a_i}b_i)_{i=1}^{n-1}$.
Thus by the induction hypothesis, $\bch {\penum m}\lambda $ has $\lambda$-sandwiching sequence $(B_{\bch {a_i}\lambda}\bch {b_i}\lambda )_{i=1}^{n-1}$.
Since we already have that $ \bch { \penum{m} }\lambda \leq \bch b\lambda$, in view of Proposition \ref{propNewNormF}, it suffices to show that
 \begin{enumerate}[label=(\alph*)]
\item \label{itNFPres2} $ B_{\bch{a }\lambda}  \bch{ b }\lambda  < B_{\bch{a + 1} \lambda } \bch{\penum{m}}\lambda $, and
\item \label{itNFPres1}  $ B_{\bch{a }\lambda}  \bch { b  }\lambda   < B_{\bch {\penu{a}}\lambda }( \bch { \penu{b} } \lambda + 1) $ if $n>1$.
 \end{enumerate}
From this and Proposition \ref{propNormForm} we may conclude that $(B_{\bch {a_i} \lambda } \bch {b_i} \lambda )_{i=1}^{n}$ is the $\lambda$-sandwiching sequence for $\bch m\lambda$.
 Consider two sub-cases.
 \medskip

\begin{enumerate}[label*= .\arabic*,wide, labelwidth=!, labelindent=0pt]
\item ($n>1$ and $a  + 1 = \penu{a}$).
In this case we have the inequality
\begin{equation}\label{eqStep}
\bch b\lambda  < \bch { A_{a }^{k-1} A_{\penu{a}}\penu{b}}\lambda = B_{\bch{a }\lambda}^{k-1} B_{\bch{\penu{a}}\lambda}  \bch{\penu{b}} \lambda,
\end{equation}
where the inequality is by Lemma \ref{lemmSandwiching}.\ref{itSandwichingTwo} and Proposition \ref{propMonotonicity} and the equality holds because $A_{a }^{x} A_{\penu{a}}\penu{b} $ is in normal form for all $x<k$, given Proposition \ref{propNormForm}.
Hence
$B_{\bch { a } \lambda + 1}  \bch {\penum{m} } \lambda   = B_{\bch { a } \lambda }^\lambda B_{\bch { a} \lambda + 1} \big ( B_{\bch { \penu{a} } \lambda }  \bch {\penu{b}} \lambda   - 1 \big )
 > B_{\bch { a } \lambda }^\lambda  B_{\bch { \penu{a} } \lambda } \bch {\penu{b}} \lambda 
>  B_{\bch{a }\lambda} \bch { b }\lambda$,
where the last inequality uses \eqref{eqStep}, thus establishing \ref{itNFPres2}.
For \ref{itNFPres1},
\begin{align*}
B_{\bch{a }\lambda}  \bch { b  }\lambda 
& <  B_{\bch{a }\lambda} B_{\bch{a }\lambda}^{k-1} B_{\bch{\penu{a}}\lambda}  \bch{\penu{b}} \lambda \\
& < B_{\bch{\penu{a}}\lambda -1}^\lambda B_{\bch{\penu{a}}\lambda}  \bch{\penu{b} } \lambda 
  =  B_{\bch{\penu{a}}\lambda} (\bch{\penu{b} } \lambda + 1),
\end{align*}
where the first inequality uses \eqref{eqStep} and the second $k<\lambda$ and $\bch a \lambda \leq \bch {\penu a} \lambda - 1$.
  \medskip

\item ($n=1$ or $a  + 1 < \penu{a}$).
Let us establish \ref{itNFPres2} first.
The argument is similar to that for \ref{itCritPropMon} of Proposition \ref{propMonotonicity}.  
Let $t \in [1 , \penum m + 1] $ be least with the property that
$A_{a + 1} (\penum{m}-t) < A_a \penum{m} ,$
and similarly $r < k$ be greatest such that
$u  := A^r_{a } A_{a + 1} (\penum{m}  - t ) < A_a \penum{m}  $.
Since $a+1< \penu a$, $A_{a} \penum m$ and $A_{a+1} \penum m$ are in normal form by Corollary \ref{corSmallIndex}.
By Lemma \ref{lemmNormFormZero} we have that $A_{a+1} \penum m = A_a^{tk} A_{a+1} (b - t) $, hence if $1 \leq s < t k - r$,
$ A_a \penum m \leq A_a^s u < A_{a+1} \penum m $,
so that by Lemma \ref{lemmNormForm}.\ref{itSmall}, $A_a^s u $ is in normal form.
Since
$A_{a } b  < A_{a + 1} \penum{m} =  A^{t k - r} _a u $, 
we have that
$b < A_{a }^{tk-r - 1}  u.$
Then,
\begin{align*}
 B_{\bch {a } \lambda} &  \bch {b } \lambda 
  < B_{\bch {a } \lambda}  \bch{  A_{a }^{tk-r -1}  u } \lambda 
 = B^{t k - r }_{\bch { a } \lambda }
  \bch{  u } \lambda  
   < B^{t k  }_{\bch { a } \lambda } \bch{ A_{a } \penum{m} } \lambda  
  = B^{t k + 1 }_{\bch { a } \lambda }  \bch{ \penum{m}} \lambda   \\
  &    < B^{t k  + 1 }_{\bch { a } \lambda }
 B^{t-1 }_{\bch { a } \lambda } B_{\bch { a } \lambda + 1  } ( \bch{ \penum{m} } \lambda - t  )  
   = B^{t \lambda }_{\bch { a } \lambda }   B_{\bch { a } \lambda + 1 } ( \bch{ \penum{m}} \lambda - t  ) 
    = B_{\bch { a } \lambda + 1 }  \bch{ \penum{m}} \lambda  .
\end{align*}
In the case that $n>0$, \ref{itNFPres1} immediately follows, since
$B_{\bch{\penu{a}}\lambda} (\bch {\penu{b}}\lambda + 1)  = B_{\bch{\penu{a}}\lambda - 1}^\lambda B_{\bch{\penu{a}}\lambda}  \bch {\penu{b}}\lambda 
> 
B_{\bch { a } \lambda + 1 }  \bch{ \penum{m}} \lambda  
>  B_{\bch {a } \lambda}   \bch {b } \lambda  $.\qedhere
\end{enumerate}
\end{enumerate}
\endproof

\section{Proofs of the Main Theorems}

In this section we put our results together to prove Theorems \ref{theoMain} and \ref{theoIndep}.
The next lemma is the final ingredient we need for the former.
\longversion{
It states that base change to $\omega$ gives the same result before and after finite base change.
}

\begin{lemma}\label{lemmBchPsiCommute}
If $k<\ell <\omega$ and $m<\omega$ then $\pk m= \ps {\ell} \bch{m}{\bases k\ell}$.
\end{lemma}

\proof
By induction on $m$.
The claim is trivial if $m=0$, so we assume otherwise. Let $m\gknf A_ab + c$ and write $\bch x\ell $ for $\bch x {\bases k\ell}$. Then $\bch m{\ell}$ has the $\ell$-normal form $ B_{\bch{a}{\ell}}\bch b{\ell} + \bch c \ell$, and the induction hypothesis yields
$    \ps {\ell}\bch m { \ell}
 =   \ot _{\ps {\ell}\bch a {\ell}} \ps {\ell}\bch b \ell +  \ps \ell \bch c \ell
 =^{\text{\sc ih}}   \ot _{\ps {k}a} \ps {k} b + \pk c =   \ps {k} m$.
\endproof

With this we are ready to prove termination for the Ackermannian Goodstein process, based on the fact that the sequence $\ok mi$, as defined below, is a decreasing sequence of ordinals.

\begin{definition}
Given $m,i\in \mathbb N$, we define $\ok mi = \ps {i+2} \good mi$.
\end{definition}

\proof[Proof of Theorem \ref{theoMain}]
Let $(\good mi)_{i<\al}$ be the Ackermannian Goodstein sequence starting on $m$ and $i$ be such that $i+1< \al$. Then,
\begin{align}
\nonumber \ok m{i+1} & = \ps {i+3} \good m{i+1}  = \ps {i+3} ( \bch {\good m {i}}{i+3} - 1) \\
\nonumber &< \ps {i+3}  \bch {\good m {i}}{i+3}  = \ps { i+2} \good m i 
 = \ok mi,
\end{align}
where the inequality follows from Proposition \ref{propMonotonicity} and the second-to-last equality from Lemma \ref{lemmBchPsiCommute}.
Hence $(\ok mi)_{i < \al}$ is decreasing below $\varepsilon_0$, so $\al$ must be finite.
\endproof

In Appendix \ref{secGamma} we review the relations $\peq_k$, fundamental sequences for $\Gamma_0$, and Theorem \ref{theoATRInd}, which is unprovable in $\atr$.
We show that Theorem \ref{theoMain} is also unprovable by deriving Theorem \ref{theoATRInd} from it, using the following key lemma.

\begin{lemma}\label{lemmBCHvsFS}
If $0<m<\omega$ and $1<k<\omega$ then
\[ \fs{\pk m }k \peq _1 \pskpe(\bch m{k+1}-1).\]
\end{lemma}

\proof
Let $\ell = k+1$ and proceed by induction on $m$.
We write $A_a b $ for $A_a(k,b)$
and $B_a b $ for $A_a(\ell,b)$.
It suffices to show
$\fs{\pk m}k \leq \ps\ell(\bch m{\ell}-1)$.
The refined inequality follows from the Bachmann property (Proposition \ref{propBachmann}),
since Proposition \ref{propMonotonicity} yields
$\pk m = \ps\ell \bch m\ell >\ps\ell( \bch m\ell-1)$.
Write $m \gknf A_ab + c$.
Note that by Proposition \ref{propNFP}, $\bch m\ell \equiv_\ell B_{\bch a \ell } \bch b \ell + \bch c \ell$.
We consider several cases.
\medskip

\begin{enumerate}[label*={\sc Case \arabic*},wide, labelwidth=!, labelindent=0pt]

\item ($a=b=c=0 $).
\  $\ps\ell( \bch m \ell -1)=0=\fs{\phi_00}k=\fs {\pk m} k $.
\medskip

\item ($c>0$). \  $\ps  \ell  ( \bch m\ell-1)
=\ps\ell (B_{\bch a\ell} \bch b\ell +\bch c \ell -1)$
\begin{align*}
&=  \ot_{\ps\ell \bch a\ell }\ps\ell \bch b\ell +\ps\ell (\bch c\ell -1)
\geq^{\text{{\sc ih}}}  \ot_{\pk a }\pk b   +\fs{\pk  c }k\\
&=  \fs{\big (\ot_{\pk a }\pk b   +\pk  c \big )} k
= \fs{ \pk m}k
.
\end{align*}

\item ($a = c = 0$ and $b> \penum{m}$).
Lemma \ref{lemmNormBCases} and Proposition \ref{propMonotonicity} imply that $\pk b\not\in \Fix0$, since $b = A_de +s$ with either $d \leq 0 $ or $s>0$, which means that $\pk b = \ot_{\pk d} \pk e + \pk s$ with $\pk d = 0$ or $ \pk b > \pk s>0$.
Therefore, $\fs{\ot_0\pk b}k=\ot_0 \fs {\pk b}k $. By Lemma \ref{lemmNormFormMinusOne}.\ref{itBMinusOne}, $B_0 (\bch b\ell -1 )$ is in normal form, hence by Corollary \ref{corMultK}, $B_0 (\bch b\ell -1 ) \cdot k$ is in extended $\ell$-normal form, so that
\begin{align*}
\ps \ell (\bch m\ell -1) & =\ps \ell ( \ell^{\bch b\ell}-1)
>\ps \ell (  \ell ^{ \bch b\ell -1} \cdot k )\\
&= \big (  \phi_0 \ps \ell (\bch b \ell -1) \big ) \cdot k \geq^{\text{\sc ih}} ( \phi_0\fs {\pk b}k )\cdot k 
> \fs {\pk m} k .
\end{align*}

\item ($a = c = 0$ and $b = \penum{m} > 0$).
Lemma \ref{lemmNormBCases} yields $b \gknf A_de$ for some $d>a$, so that $\pk b\in \Fix0$ and $ \fs { \pk m}k =  \pk b   \cdot k 
$.
By Corollary \ref{corMultK}, $ \bch b \ell \cdot k$ is in simplified $\ell$-normal form and hence
$\ps \ell(\bch m \ell-1)  = \ps \ell (\ell^{\bch b \ell}-1 ) 
 > \ps \ell({\bch b \ell} \cdot k)
=  \pk b   \cdot k  
= \fs{\pk m}k$,
where the inequality follows from $(x+1)^y > xy 
$.
\medskip

\item ($a>0$ and $b=c=0$).
By Lemma \ref{lemmNFLeftExp}, $B_{ \bch a\ell-1}^{j} 1$ is in $\ell$-normal form for $1\leq j\leq k$.
Then Proposition \ref{propMonotonicity} and the induction hypothesis yield 
\begin{align*}
&\ps \ell ( \bch m \ell -1)
=\ps \ell(B_{\bch a \ell } 0-1)
=\ps \ell(B_{ \bch a \ell -1}^{\ell} 1-1)\\
&> \ps \ell B_{ \bch a \ell -1}^{k} 1  
 =  \ot_{ \ps \ell ( \bch a \ell -1)}^k 1
 \geq^{\text{\sc ih}} \ot_{\fs {\pk a} k }^k 1 
 \geq  \fs{\ot_{\pk a }0}k = \fs{\pk m}k.
\end{align*}

\item ($a>0$, $b>\penum{m}$ and $c= 0$).
Once again by Proposition \ref{propMonotonicity} and Lemma \ref{lemmNormBCases} one can check that $\pk b\not\in \Fix{\pk a}$ and thus
\begin{align*}
&\fs {\ot_{\pk a  }\pk b}k   =  \ot_{\pk a } \fs{\pk b} k  
\leq^{\text{\sc ih}}   \ot_{\ps \ell \bch a \ell } {\ps \ell ( \bch b \ell - 1) }\\
& = \ps \ell  B_{\bch a\ell } ( \bch b\ell - 1)
 < \ps \ell ( B_{\bch a\ell }  \bch b\ell - 1)
  = \ps \ell ( \bch m \ell - 1).
\end{align*}

\item ($a>0$, $b= \penum{m}$ and $c = 0$).
Note that in this case $\pk b\in \Fix{\pk a}$ so that $\fs{\ot_{\pk a }\pk b}k = \ot_{\fs{\pk a }k}^{\xcases k{\pk a}} \pk b .$
Moreover, note that $B_{\bch a \ell -1}^{j}\bch b \ell $ is in $ \ell $-normal form for $j\leq k$ by Corollary \ref{corSmallIndex}.
Then Proposition \ref{propMonotonicity} and  the induction hypothesis yield
\begin{align*}
\ps \ell(\bch m \ell-1)
&=\ps \ell(B_{\bch a \ell} \bch b \ell -1)
=\ps \ell \big (B_{\bch a \ell-1}^{k+1}  B_{\bch a \ell}(\bch b \ell-1)-1 \big )\\
& \geq \ps \ell B_{\bch a \ell-1}^{k} \bch b \ell
=  \ot_{\ps \ell(\bch a \ell-1)}^k \ps \ell \bch b \ell\\
&\geq^{\text{\sc ih}} \ot_{ \fs {\pk a}k }^k \pk b
\geq  \fs{ \ot_{\pk a } \pk b}k.&\qedhere
\end{align*}

\end{enumerate}

\longversion{Since we have considered all possible cases, the claim follows.}
\endproof

With this we may show that $\ok mk$ decreases at least as slowly as $\fsi{\ok m0}{k+1}$.

\begin{lemma}\label{lemmOvsFS}
If $ m\geq 0$ and $k\geq 2$ then $\ok mk\seq_{k+1} \fsi{\ok m0}{k+1}$.
\end{lemma}

\proof\label{lemmODecr}
The assertion holds for $k=0$.
Assume now that
$\ok mk\seq_{k+1} \fsi {\ok m0}{k+1}.$
Then $\ok mk\seq_{k+2} \fsi{\ok m0}{k+1}$  and hence
$\fs{\ok mk}{k+2}\seq_{k+2} \fs { \fsi {\ok m0}{k+1}} {k+2}$.
Moreover Lemma \ref{lemmBCHvsFS} yields
$\ok m{k+1} 
= \ps {k+3} \good m{k+1} 
 = \ps {k+3} (\bch{ \good m{k }  }{k+3}-1) 
  \seq_1 \fs{ \ps { k+2 } \good m{k }  }{k+2}
= \fs{\ok mk}{k+2}$,
hence 
$\ok m{k+1} \seq_{k+2} \fs{\ok mk}{k+2}$.
Putting things together we arrive at
$\ok m{k+1} \seq_{k+2}  \fsi{\ok m0}{k+2}.$
\endproof

\proof[Proof of Theorem \ref{theoIndep}]
Let $\gamma_n$ be as in Theorem \ref{theoATRInd}.
Moreover let $a_0:=0$ and $a_{n+1}:=A_{a_n}(2,0).$ Recall that $\ok {a_m} \ell =\ps  {\ell +2}\good {a_m}\ell $.
Then $\ok {a_m}0 =\ps  2 a_m =\gamma_m $ for $m\geq 3$.
Lemma 10 yields 
$\fsi {\ok {a_m}0}{k+1}\peq_{k+1} \ok{a_m}k $.

Assume that
$\atr \vdash \forall m \exists \ell \ (\good m \ell =0 ) .$
The function $m\mapsto a_m$ is provably computable in $ \atr$ and therefore
$\atr \vdash \forall m \exists \ell \ ( \good {a_m}\ell =0 ).$
From this we obtain
$\atr\vdash \forall m \exists \ell \ ( \ok{a_m}\ell=0 ).$
But then by formalizing the proof of Lemma \ref{lemmOvsFS}, we see that
$\atr \vdash \forall m \exists \ell \ ( \fsi{\ga_m}{\ell}=0 ) ,$
contradicting Theorem \ref{theoATRInd}.
\endproof

\appendix

\section{The Feferman-Sch\"utte Ordinal}\label{secGamma} 

\longversion{
We will assume some familiarity with notation systems for $\Gamma_0$; the interested reader may find a more detailed treatment in e.g.~\cite{Pohlers:2009:PTBook}.
We will use the modified fixpoint-free Veblen hierarchy, but it is convenient to introduce the standard hierarchy first.

Given an ordinal $\beta$ we define a function $\varphi_\beta$ recursively as follows.
Let $\varphi_0\al:=\om^\al$ and for $\be>0$ let $\varphi_\be\ga$ be the $\ga$-th member of
$\{\eta:(\forall\xi<\be) [\varphi_\xi \eta=\eta]\}$.}
\shortversion{
In this appendix we briefly review the ordinal $\Gamma_0$. For a more detailed treatment see e.g.~\cite{Pohlers:2009:PTBook}.
Given an ordinal $\alpha$ recall that $\varphi_\alpha$ is defined recursively so that $\varphi_0\beta:=\om^\beta$ and for $\alpha>0$, $\varphi_\al\be$ is the $\beta$-th member of
$\{\eta:(\forall\xi<\be) [\varphi_\xi \eta=\eta]\}$.}
We will use the fixed point-free variant given by
\[\ot_\al\be:=
\begin{cases}
\varphi_\al(\be+1)&\mbox{ if there exists a $\beta_0$ such that $\be=\be_0+n$ and $\varphi_\al\be_0=\be_0$},\\
\varphi_\al\be&\mbox{ otherwise.}
\end{cases}\]
Then, $\Gamma_0$ is the first non-zero ordinal closed under $(\alpha,\beta)\mapsto \ot_\alpha\beta$.\david{I thought this notation looked nicer than $\overline\varphi$, but we can change it back if it is too non-standard.}
The technical benefit of the modified hierarchy is witnessed by the following.
\begin{proposition}\label{propGammaNF}
If $0<\xi<\Gamma_0$, there exist unique $\alpha,\beta,\gamma<\xi$ such that $\xi = \ot_\alpha\beta+\gamma$.
\end{proposition}

We will call this the {\em Veblen normal form} of $\xi$ and write $\xi \vnf \ot_\alpha \beta+\gamma$.
Proposition \ref{propGammaNF} does not hold for $\varphi$, as for example $\varphi_0\varphi_1 0 = \varphi_1 0 < \Gamma_0$.
The order relation between elements of $\Gamma_0$ can be computed recursively on their Veblen normal form.
Below we consider only $\xi,\zeta>0$, as clearly $0 < \ot_\alpha\beta + \gamma$ regardless of $\alpha,\beta,\gamma$.

\begin{lemma}
Given $\xi,\xi'<\Gamma_0$ with $\xi=\ot_{\alpha}{\beta }+{\gamma }$ and $\xi' = \ot_{\alpha'}{\beta'}+{\gamma'}$ both in Veblen normal form, $\xi<\xi'$ if and only if
\begin{multicols}2
\begin{enumerate}

\item $\alpha=\alpha'$, $\beta=\beta'$ and $\gamma<\gamma'$;

\item $\al<\al'$ and $\be<\ot_ {\alpha'}\beta'$;

\item $\al=\al'$ and $\be<\be'$, or

\item $\al'<\al$ and $\ot _\al\be\leq \be'$.

\end{enumerate}
\end{multicols}
\end{lemma}

The above properties of modified Veblen functions will suffice to prove that the Ackermannian Goodstein process always terminates on finite time.
In order to prove independence from $\atr$, we also need to review fundamental sequences.
Let
$\Fix\al:=\{\ot_\be\ga:\ga>\al\}$; these are the fixed points of $\varphi_\al$. Fundamental sequences will be defined separately at such points.

\begin{definition}
For $\xi < \Gamma_0$ and $x<\omega$, define $\fs \al x$ recursively as follows.
First we set $ \xcases x \alpha = x$ if $\alpha$ is a successor, $ \xcases x\alpha = 1$ otherwise.
Then, define:
\begin{multicols}2
\begin{enumerate}

\item $\fs 0 x = 0$.

\item $\fs{ ( \ot_\al \be + \gamma ) }x= \ot_\al \be +  \fs \ga x$ if $\ga > 0$.

\item $\fs{\phi_00} x =0$ (note that $  {\phi_00} = 1$).


\item $\fs{\phi_0\la} x =\la\cdot x$ if $\la \in \Fix0$.

\item $\fs{\ot_\al0}x:=\ot^{\xcases x\alpha}_{\fs \al x }  1 $ if $\al > 0$.

\item $\fs{\ot_\al(\be+1)}x:=
\ot^{\xcases x\alpha}_{\fs\al x} \ot_\al\be $ if $\al > 0 $.

\item $\fs{\ot_\al\la}x:=\ot_{\al}\fs \la x$ if $\la$ is a limit and $\la\not\in \Fix\al$.

\item $\fs{\ot_\al\la}x:=\ot^{\xcases x\alpha}_{\fs \al x} \la $ if $\la\in \Fix\al$ and $\al>0$.
 
\end{enumerate}
\end{multicols}

\end{definition}

For an ordinal $\xi<\Gamma_0$ and $n\geq 2$ we define inductively $\fsi \xi 2 = \fs \xi 2$ and $\fsi\xi{ n+1} = \fs{\fsi \xi  n }{n+1}$. Define $\prec_k$ to be the transitive closure of $\{(\al,\be):\al=\fs \be k\}$, and $\peq_k$ to be its reflexive closure.
The ordering $\peq_1$ satisfies the {\em Bachmann property:}

\begin{proposition}\label{propBachmann}
If $\alpha,\beta < \Gamma_0$ and $k< \omega$ are such that $\fs \al k < \beta < \al$, then $\fs\al k\peq_1 \be$.
\end{proposition}

See, for example, \cite{BCW}, \cite{Weiermanna} and \cite{WeiermannLC} for more details.
It is further well-known that $\al\peq_k\be $ yields $\al\peq_{k+1}\be$.
Our unprovability result for $\atr$ follows from Theorem \ref{theoATRInd} below; see any of \cite{Buchholza,FairtloughWainer,WeiermannBSL} for a proof.

\begin{theorem}\label{theoATRInd}
By recursion on $n$ define $\gamma_n<\Gamma_0$ given by $\gamma_0:=0$ and $\gamma_{n+1}:=\ot_{\ga_n} 0 .$
Then, $ \forall m \exists \ell \ ( \fsi{\ga_m}\ell  =0 )$ is not provable in $ \atr $.
\end{theorem}

\section{Alternative Normal Forms}\label{secConc}

There are other notions of normal form that we may consider aside from those of Definition \ref{defNF}.
For example, we may skip the sandwiching procedure and choose $a_0,b_0$ so that $A_{a_0}{b_0}$ is maximal with the property that $A_{a_0}b_0\leq m$, then choose $a$ maximal such that there exists $b\geq 0$ with $A_ab =A_{a_0}{b_0}$; call this the {\em alternative $k$-normal form} of $m$.

These alternative normal forms can be rather inefficient.
Fix a base $k$. We claim that for any $j > 0$ and any $a,b$, if $A_ab = A^j_1 2$ then $a\leq 1$.
It suffices to show that $A^j_1 2$ is not in the range of $A_2$.
Proceed by induction on $j$, and towards a contradiction, assume that $A^j_1 2 = A_2 d$.
Then $A^j_1 2 = A^k_1 A_2 (d-1)$.
Now consider two cases.
If $j\leq k$ then $ 2 = A^{k-j}_1 A_2 (d-1)$.
If $d=k-j=0$ then
$ A^{k-j}_1 A_2 (d-1) =A_2 (-1) = 1 < 2$, 
otherwise
$A^{k-j}_1 A_2 (d-1) \geq A _1 1 > 2$.
If $j>k$ then $A^{j-k}_1 2 = A_2 (d-1)$, contradicting the induction hypothesis if $d>0$ and false if $d=0$ since $A_2(-1) = 1$.
Thus under such normal forms $m = A^j_1 2$ would have to be written with $j$ instances of $A_1$, even if $m$ is larger than (say) $A_{100} 0$.
Our sandwiching procedure is tailored to avoid such a situation.

A third approach would be to choose the normal form of $m$ to be the shortest possible with respect of the number of symbols used. However, we currently do not know if given an expression $\tau$, there is a primitive recursive procedure to compute the shortest $\tau_\ast$ with the same value as $\tau$.

If $\tau$ is a formal expression built from $0$, $A_xy$, and $+$, define $\|\tau\|$ by $\|0\|= 1$ and $\|A_\pi\sigma  + \rho\| = \|\pi\|+\|\sigma\|+\|\rho\|$.
Neither the original nor the alternative $k$-normal forms produce minimal norms.
Let $m =A_0 A_10 -1$.
Then $m  \enfpar 2 A_10 \cdot p + q$ with a large $p $ of about $ k^{A_10 - \log_k A_10}$. 
However, $ m = \sum_{i=0}^{A_10-1} A_0 i $ which has norm about $(A_10)^2 $.
Now let $n=A_1(A_1^{k-1}A_2 (k-1) + 1 )$. Then $n$ is in alternative normal form with norm about $4k$, since $k-1$ can only be written as $A_00\cdot (k-1)$.
But $n=A_0^k A_2 A_01 $ which has norm about $k$.

Nevertheless, Goodstein sequences for such normal forms could be of interest and we leave their study as an open line of research.
We conjecture that Ackermannian Goodstein sequences will terminate for any `reasonable' notion of normal form.

\bibliographystyle{amsplain}

\bibliography{biblio}

\providecommand{\bysame}{\leavevmode\hbox to3em{\hrulefill}\thinspace}
\providecommand{\MR}{\relax\ifhmode\unskip\space\fi MR }
\providecommand{\MRhref}[2]{%
  \href{http://www.ams.org/mathscinet-getitem?mr=#1}{#2}
}
\providecommand{\href}[2]{#2}
\begin{thebibliography}{10}

\bibitem{Buchholza}
W.~Buchholz and S.~Wainer, \emph{Provably computable functions and the fast
  growing hierarchy}, Contemporary Mathematics: Logic and combinatorics
  (Providence, RI), vol.~65, American Mathematical Society, 1987, pp.~179--198.

\bibitem{Cichon}
E.A. Cichon, \emph{A short proof of two recently discovered independence
  results using recursion theoretic methods}, Proceedings of the American
  Mathematical Society \textbf{87} (1983), 704--706.

\bibitem{BCW}
E.A. Cichon, W.~Buchholz, and A.~Weiermann, \emph{A uniform approach to
  fundamental sequences and hierarchies}, Mathematical Logic Quarterly
  \textbf{40} (1994), 273--286.

\bibitem{FairtloughWainer}
M.~Fairtlough and S.S. Wainer, \emph{Hierarchies of provably recursive
  functions}, Handbook of Proof Theory (Samuel~R. Buss, ed.), Elsevier Science
  BV, 1998, pp.~149--207.

\bibitem{Feferman:1964:SystemsOfPredicativeAnalysis}
S.~Feferman, \emph{Systems of predicative analysis}, Journal of Symbolic Logic
  \textbf{29} (1964), 1--30.

\bibitem{Feferman}
\bysame, \emph{Systems of predicative analysis, {II:} representations of
  ordinals}, Journal of Symbolic Logic \textbf{33} (1968), no.~2, 193--220.

\bibitem{Pelupessy}
H.~Friedman and F.~Pelupessy, \emph{Independence of {R}amsey theorem variants
  using $\varepsilon_0$}, Proceedings of the American Mathematical Society
  \textbf{144} (2016), no.~2, 853--860, English summary.

\bibitem{Gentzen1936}
G.~Gentzen, \emph{Die {W}iderspruchsfreiheit der reinen {Z}ahlentheorie},
  Mathematische Annalen \textbf{112} (1936), 493--565.

\bibitem{Godel1931}
K.~G\"{o}del, \emph{{\"{U}ber Formal Unentscheidbare S\"{a}tze der Principia
  Mathematica und Verwandter Systeme, I}}, Monatshefte f\"{u}r Mathematik und
  Physik \textbf{38} (1931), 173--198.

\bibitem{Goodstein1944}
R.L. Goodstein, \emph{On the restricted ordinal theorem}, Journal of Symbolic
  Logic \textbf{9} (1944), no.~2, 33–41.

\bibitem{Goodsteinb}
\bysame, \emph{Transfinite ordinals in recursive number theory}, Journal of
  Symbolic Logic \textbf{12} (1947), no.~4, 123--129.

\bibitem{graham1969}
R.~L. Graham and B.~L. Rothschild, \emph{Ramsey's theorem for $n$-dimensional
  arrays}, Bulletin of the American Mathematical Society \textbf{75} (1969),
  no.~2, 418--422.

\bibitem{Kirby}
L.~Kirby and J.~Paris, \emph{Accessible independence results for {P}eano
  arithmetic}, Bulletin of the London Mathematical Society \textbf{14} (1982),
  no.~4, 285--293.

\bibitem{ParisHarrington}
J.~Paris and L.~Harrington, \emph{A mathematical incompletenss in {P}eano
  arithmetic}, Handbook of Mathematical Logic (J.~Barwise, ed.), North-Holland
  Publishing Company, 1977, pp.~1133--1142.

\bibitem{Pohlers:2009:PTBook}
W.~Pohlers, \emph{Proof theory, the first step into impredicativity},
  Springer-Verlag, Berlin Heidelberg, 2009.

\bibitem{Simpson:1985:FriedmansResearh}
S.G. Simpson, \emph{Friedman's research on subsystems of second order
  arithmetic}, Harvey Friedman's Research in the Foundations of Mathematics
  (L.~Harrington, M.~Morley, A~\v{S}\v{c}edrov, and S.~G. Simpson, eds.),
  North-Holland, 1985, pp.~137--159.

\bibitem{Simpson:2009:SubsystemsOfSecondOrderArithmetic}
\bysame, \emph{Subsystems of {S}econd {O}rder {A}rithmetic}, Cambridge
  University Press, New York, 2009.

\bibitem{Weiermanna}
A.~Weiermann, \emph{Some interesting connections between the slow growing
  hierarchy and the ackermann function}, Journal of Symbolic Logic \textbf{66}
  (2001), no.~2, 609--628.

\bibitem{WeiermannLC}
\bysame, \emph{A very slow growing hierarchy for $\gamma_0$}, Logic Colloquium
  '99 (Urbana, IL), Lecture Notes in Logic, vol.~17, 2004, pp.~182--199.

\bibitem{WeiermannBSL}
\bysame, \emph{Classifying the provably total functions of {PA}}, Bulletin of
  Symbolic Logic \textbf{12} (2006), no.~2, 177--190.

\bibitem{vanHoof}
A.~Weiermann and W.~Van~Hoof, \emph{Sharp phase transition thresholds for the
  {P}aris {H}arrington {R}amsey numbers for a fixed dimension}, Proceedings of
  the American Mathematical Society \textbf{140} (2012), no.~8, 2913--2927.

\bibitem{Weyl}
H.~Weyl, \emph{The continuum: A critical examination of the foundation of
  analysis}, Mineola: Dover, 1918.

\end{thebibliography}

\end{document}